\documentclass[preprint,12pt]{elsarticle}

\usepackage{amssymb}
\usepackage{mathrsfs}
\usepackage{stmaryrd}
\input diagxy

\newtheorem{Th}{\bf Theorem}[section]
\newtheorem{Lem}[Th]{\bf Lemma}
\newtheorem{Pro}[Th]{\bf Proposition}

\newtheorem{Def}[Th]{\bf Definition}
\newtheorem{Exam}[Th]{\bf Example}
\newtheorem{Rem}[Th]{\bf Remark}

\usepackage{amssymb}

\newcommand{\p}[1]{{\bf Proof.} #1 \ $\Box$}

\journal{Advances in Mathematics}

\begin{document}

\begin{frontmatter}



\title{Algebraic representation of $L$-valued continuous lattices via the open filter monad}


\author{Wei Yao}

\address{School of Mathematics and Statistics, Nanjing University of Information Science and Technology, Nanjing, China}

\author{Yueli Yue}

\address{School of Mathematical Sciences, Ocean University of China, Qingdao, China}

\author{Bin Pang}

\address{School of Mathematics and Statistics, Beijing Institute of Technology, Beijing, China}

\begin{abstract}
With a complete Heyting algebra $L$ as the truth value table, we prove that the collections of open filters of stratified $L$-valued topological spaces form a monad. By means of $L$-Scott topology and the specialization $L$-order, we get that the algebras of open filter monad are precisely $L$-continuous lattices.

\end{abstract}


\end{frontmatter}

\section{Introduction}

Monad (also called triple) is an endofunctor over a base category, together with two natural transformations
satisfying two coherence conditions \cite{Monad}. The categories of algebras arising from a monad
seem to be the most natural generalization of Birkhoff's equational classes \cite{Birkhoff}.
Nowadays, the monad approach becomes a very important way to study topological structures.

In \cite{Man}, Manes proved that compact Hausdorff spaces can be obtained as the Eilenberg-Moore algebras of the ultrafilter monad. Later Barr \cite{Bar} showed that by going from {\bf Set} (the category of sets) to {\bf Rel} (the category of binary relational sets) and relaxing the axioms on the monad, and then the related algebras derived, it was possible to obtain all topological spaces as lax algebras of a suitable lax extension of the ultrafilter monad to {\bf Rel}.

Domain Theory initiated by Dana Scott in the late seventies, from a viewpoint of pure mathematics,
can be considered as a combined branch of order/lattice theory, topology, logic, category theory and so on.
In \cite{Scott}, Dana Scott characterized the continuous lattices endowed with the Scott topology precisely as the spaces that are injective over all subspace embedding. This was the first important result of connections between algebraic structures and topological structure in domain theory in the literature.
For another perspective, Alan Day \cite{Day} and Oswald Wyler \cite{Wyler} independently characterized the continuous lattices as the algebras of the open filter monad on the category of $T_0$ topological spaces.
One uses the fact that the filter monad is of Kock-Z\"{o}berlein type, and that in any poset-enriched category with such a monad structure, the injective objects over a certain class of embeddings defined in terms of monad structures are precisely the algebras. The conclusion follows from the fact that the embeddings associated to the filter monad are exactly topological embeddings.
Escardo and Flagg \cite{Escardo} studied various kinds of injective spaces over different kinds of topological embeddings.

Quantitative domains are extensions of classical domains by generalizing ordered relations to other more general structures, such as enriched categories, generalized metrics and lattice-valued fuzzy ordered relations.
Lattice-valued fuzzy order approach to domain theory is originated by Fan and Zhang \cite{FanTCS,FanZhang}, and then mainly developed by Yao \cite{PartI,PartII},
Zhang \cite{Zhang1,Zhang2,Zhang4,Zhang3,Zhang5}, Li Qingguo \cite{LiQG5,LiQG2,LiQG1,LiQG3,LiQG4}, Zhao \cite{Zhao2,Zhao3,Zhao1}, etc.

In \cite{YaoTFS}, Yao showed that frame-valued continuous lattices are categorical isomorphic to injective lattice-valued $T_0$ spaces. The aim of this paper is to generalize Day's result to lattice-valued setting to show that frame-valued continuous lattices are exactly the algebras of the open filter monad on the category of $T_0$ frame-valued topological spaces.

\section{Preliminaries}

\subsection{The truth value table}

We refer to \cite{Domain} for the contents of lattice theory.

A poset $L$ is called a {\it complete lattice} if the supermum $\bigvee S$ exists for each $S\subseteq L$, or equivalently, the infimum $\bigwedge T$ exists for each $T\subseteq L$. For the case of $\emptyset$, $\bigvee\emptyset$ and $\bigwedge\emptyset$ are the least and the greatest elements, denoted by $0$ and $1$, respectively.

\begin{Def}A complete lattice $L$ is called a frame, or a complete Heyting algebra, if it satisfies the first infinite distributive law, that is,
$$a\wedge\bigvee S=\bigvee_{s\in S}a\wedge s\ (\forall a\in L,\ \forall S\subseteq L).$$\end{Def}

For a frame $L$, the related implication operation $\rightarrow:L\times L\longrightarrow L$ is given by $$a\rightarrow b=\bigvee\{c\in L|\ a\wedge c\leq b\}\ (\forall a,b\in L)$$Then we get an adjoint pair $(\wedge,\rightarrow)$ satisfying that $$c\leq a\rightarrow b\Longleftrightarrow a\wedge c\leq b\ (\forall a,b,c\in L).$$

In this paper, $L$ always denotes a frame.

\begin{Pro}For all $a,b,c\in L,\ \{a_i|\ i\in I\},\ \{c_i|\ i\in I\},\ \{b_j|\ j\in J\}\subseteq L$, it holds that

{\rm (1)} $a\rightarrow b=1\Longleftrightarrow a\leq
b$;

{\rm (2)} $1\rightarrow a=a$;

{\rm (3)} $a\wedge(a\rightarrow b)=a\wedge b$;

{\rm (4)} $b\leq a\rightarrow (a\wedge b)$, $a\leq (a\rightarrow b)\rightarrow b$;

{\rm (5)} $(a\rightarrow b)\wedge (b\rightarrow c)\leq a\rightarrow c$;

{\rm (6)} $(\bigvee_ia_i)\rightarrow b=\bigwedge_i(a_i\rightarrow
b)$;

{\rm (7)} $a\rightarrow(\bigwedge_jb_j)=\bigwedge_j(a\rightarrow
b_j)$;

{\rm (8)} $(\bigwedge_ia_i)\rightarrow(\bigwedge_ic_i)\geq\bigwedge_i(a_i\rightarrow
c_i)$;

{\rm (9)} $(\bigvee_ia_i)\rightarrow(\bigvee_ic_i)\geq\bigwedge_i(a_i\rightarrow
c_i)$;

{\rm (10)} $(c\rightarrow a)\rightarrow (c\rightarrow b)\geq
a\rightarrow b$;

{\rm (11)} $(a\rightarrow c)\rightarrow (b\rightarrow c)\geq
b\rightarrow a$;

{\rm (12)} $a\rightarrow (b\rightarrow c) = (a\wedge b)\rightarrow
c=b\rightarrow (a\wedge c)$;

{\rm(13)} $b\rightarrow c\leq(a\wedge b)\rightarrow(a\rightarrow c)$.\end{Pro}

\begin{Exam}{\rm(1)} Every finite distributive lattice is a frame.

{\rm(2)} The unit interval $[0,1]$ is a frame.

{\rm(3)} For every ordinary topological space $(X,\mathcal{T})$, the pair $(\mathcal{T},\subseteq)$ is a frame.
\end{Exam}

\subsection{$L$-subsets and $L$-valued topological spaces}

We refer to \cite{MVTop} for the contents of $L$-subsets and $L$-valued topological spaces.

Let $X$ be a set. Every mapping $A:X\longrightarrow L$ is called an $L$-{\it subset} of $X$, denoted by $A\in L^X$.
For an element $a\in L$, the notation $a_X$ denotes the constant $L$-subset of $X$ with the value $a$, that is, $a_X(x)=a\ (\forall x\in X)$.

Let $f:X\longrightarrow Y$ be a mapping between two sets. Define $f^\rightarrow_L:L^X\longrightarrow L^Y$ and $f^\leftarrow_L:L^Y\longrightarrow L^X$ respectively by $$f_L^\leftarrow(B)=B\circ f,\ \ f_L^\rightarrow(A)(y)=\bigvee\limits_{f(x)=y}A(x).$$

\begin{Def}A subfamily $\mathcal{O}(X)\subseteq L^X$ is called an $L$-valued topology if

{\rm(O1)} $A,B\in\mathcal{O}(X)$ implies $A\wedge B\in \mathcal{O}(X)$;

{\rm(O2)} $\{A_j|\ j\in J\}\subseteq\mathcal{O}(X)$ implies $\bigvee_jA_j\in\mathcal{O}(X)$;

{\rm(O3)} $a_X\in\mathcal{O}(X)$ for any $a\in L$.\\
The pair $(X,\mathcal{O}(X))$ is called an $L$-valued topological space.
\end{Def}

Let $\mathcal{O}(X)$ be an $L$-valued topology. A subfamily $\mathcal{B}\subseteq\mathcal{O}(X)$ is called a {\it base} of $X$ if for every $A\in\mathcal{O}(X)$, there exists $\{(B_j,a_j)|\ j\in J\}\subseteq\mathcal{O}(X)\times L$ such that $A=\bigvee_jB_j\wedge (a_j)_X$.
An $L$-valued topological space $X$ is called $T_0$ if $A(x)=A(y)\ (\forall A\in\mathcal{O}(X))$ implies $x=y$.

A mapping $f:(X,\mathcal{O}(X))\longrightarrow(Y,\mathcal{O}(Y))$ between two $L$-valued topological spaces is called {\it continuous} if $f_L^\leftarrow(B)\in\mathcal{O}(X)$ for each $B\in\mathcal{O}(X)$. If $Y$ has a base $\mathcal{B}$,
then it is routine to show that $f:(X,\mathcal{O}(X))\longrightarrow(Y,\mathcal{O}(Y))$ is continuous if{}f $f_L^\leftarrow(B)\in\mathcal{O}(X)$ for every $B\in\mathcal{B}$.

Let $L$-{\bf Top}$_0$ denote the category of all $T_0$ $L$-valued topological spaces and continuous mappings as morphisms.

\subsection{Category theory}

We refer to \cite{Category,Monad} for contents of category theory.

A {\it functor} $F:{\bf A}\longrightarrow{\bf B}$ between two categories is an assignment sending every {\bf A}-object $A$ to a {\bf B}-object $F(B)$ and every {\bf A}-morphism $f:A\longrightarrow A'$ to a {\bf B}-morphism $F(f):F(A)\longrightarrow F(A')$, which preserves composition and identity morphisms. For any category {\bf A}, there is an {\it identity functor} $id_{\bf A}$ sending every $f:A\longrightarrow A'$ to itself.

Let $F,G : {\bf A}\longrightarrow{\bf  B}$ be two functors. A {\it natural transformation} $\tau$ from $F$ to $G$ (denoted by $\tau: F\longrightarrow  G$) is a function that assigns to each {\bf A}-object $A$ a {\bf B}-morphism $\tau_A : F(A)\longrightarrow G(A)$ such that the following diagram commutes, that is, for each
{\bf A}-morphism $f:A\longrightarrow A'$, $G(f)\cdot\tau_A=\tau_{A'}\cdot F(f)$.

$$\bfig
\square[F(A)`F(A')`G(A)`G(A');F(f)`\tau_A`\tau_{A'}`G(f)]
\efig$$

A {\it monad} over a category {\bf X} is a triple $(T,\eta,\mu)$ consisting of
a functor $T:{\bf X}\longrightarrow{\bf X}$ and natural transformations
$\eta : id_{\bf X}\longrightarrow  T$ and $\mu : T\circ T\longrightarrow T$
such that the following diagrams commute, that is, $\mu\cdot T\mu=\mu\cdot \mu T$ and $\mu\cdot\eta T=id_T=\mu\cdot T\eta$.

$$\bfig
\square(0,0)[TTT`TT`TT`T;T\mu`\mu T`\mu`\mu]
\morphism(1200,500)[T`TT;\eta T]
\morphism(1700,500)/<-/[TT`T;T\eta]
\morphism(1700,500)|m|<0,-500>[TT`T;\mu]
\node a(1200,500)[T]
\node b(2200,500)[T]
\node c(1700,0)[T]
\arrow|b|/->/[a`c;id_T]
\arrow|b|/->/[b`c;id_T]
\efig$$

Given a monad $(T,\eta,\mu)$ over {\bf X}, a $T$-algebra (or an Eilenberg-Moore algebra) is a pair $(X,r)$, where $X$ is an {\bf X}-object and the structured morphism $r:T(X)\longrightarrow X$ satisfies, $r\cdot T(r)=r\cdot \mu_X$ and $id_X=r\cdot \eta_X$.

$$\bfig
\square(0,0)[TT(X)`T(X)`T(X)`X;T(r)`\mu_X`r`r]
\qtriangle(1200,0)/>`>`>/[X`T(X)`X;\eta_X`id_X`r]
\efig$$

\subsection{$L$-ordered sets}

We refer to \cite{Bel,PartI,PartII} for the contents of $L$-ordered sets.

\begin{Def}A mapping $e:X\times X\longrightarrow L$ is called an $L$-order if

{\rm(E1)} $e(x,x)=1$;

{\rm(E2)} $e(x,y)\wedge e(y,z)\leq e(x,z)$;

{\rm(E3)} if $e(x,y)\wedge e(y,x)=1$, then $x=y$.\\ The pair $(X,e)$ is called an $L$-ordered set.
\end{Def}

Define $A^l,\ A^u\in L^X$ respectively by $$A^l(x)=\bigwedge\limits_{y\in X}A(y)\rightarrow e(x,y)\ {\rm and}\ A^u(x)=\bigwedge\limits_{y\in X}A(y)\rightarrow e(y,x),$$ and define ${\uparrow}x$ and ${\downarrow}x$ respectively by ${\uparrow}x(y)=e(x,y)$ and ${\downarrow}x(y)=e(y,x)$. An $L$-subset $S\in L^X$ is called a {\it lower set} (resp., {\it upper set}) if $S(x)\wedge e(y,x)\leq S(y)$ (resp., $S(x)\wedge e(x,y)\leq S(y)$) for all $x,y\in X$. Clearly, $A^l$ and ${\downarrow}x$ (resp., $A^u$ and ${\uparrow}x$) are lower (resp., upper) sets for all $A\in L^X$ and $x\in X$.

\begin{Def}{\rm(1)} An element $x$ is called a supremum of $A\in L^X$, in symbols, $x=\sqcup A$, if $e(x,y)=A^u(y)\ (\forall y\in X)$.

{\rm(2)} An element $x$ is called an infimum of $A\in L^X$, in symbols, $x=\sqcap A$, if $e(y,x)=A^l(y)\ (\forall y\in X)$.\end{Def}

Clearly, $\sqcup{\downarrow}x=x=\sqcap{\uparrow}x$ for every $x\in X$. It is easy to show the following proposition.

\begin{Pro}(Yao Specialization order){\rm(1)} If $A\leq B\leq C$ and $\sqcup A=\sqcup C=x$, then $\sqcup B=x$.

{\rm(2)} $\bigwedge\limits_{z\in X}e(x,z)\rightarrow e(y,z)=e(y,x)$.

{\rm(3)} $\bigwedge\limits_{z\in X}e(z,x)\rightarrow e(z,y)=e(x,y)$.\end{Pro}

\begin{Def}An $L$-ordered set $(X,e)$ is called complete if every $L$-subset has a supremum, or equivalently, every $L$-subset has an infimum.\end{Def}

The equivalence in Definition 2.8 is because $\sqcup A=\sqcap A^u,\ \sqcap A=\sqcup A^l$ hold for $A\in L^X$ in an $L$-ordered set $X$.

\begin{Exam}
{\rm(1)} Let $(X,\mathcal{O}(X))$ be a $T_0$ $L$-valued topological space. Define $e_{\mathcal{O}(X)}:X\times X\longrightarrow L$ by $$e_{\mathcal{O}(X)}(x,y)=\bigwedge\limits_{A\in\mathcal{O}(X)}A(x)\rightarrow A(y)\ (\forall x,y\in X).$$ Then $e_{\mathcal{O}(X)}$ is an $L$-order on $X$, called the specialization $L$-order of $(X,\mathcal{O}(X))$.

{\rm(2)} Define a mapping $e_L:L\times
L\longrightarrow L$ by $e_L(x,y)=x\rightarrow y$.
Then $(L,e_L)$ is a complete $L$-ordered set, where for $A\in L^L$, $$\sqcup A=\bigvee\limits_{a\in L}a\wedge A(a),\ \sqcap A=\bigwedge\limits_{a\in L}A(a)\rightarrow a.$$Let $S,A\in L^X$ be two $L$-subsets.
By considering $S$ as a mapping from $X$ to $L$, it is easy to verify that $\sqcup S^\rightarrow_L(A)=\bigvee_{x\in X}S(x)\wedge A(x)$.

{\rm(3)} Let $X$ be a nonempty set. Define a mapping ${\rm sub}_X:L^X\times L^X\longrightarrow L$ by $${\rm sub}_X(A,B)=\bigwedge\limits_{x\in X}A(x)\rightarrow B(x).$$ Then
$(L^X,{\rm sub}_X)$ is a complete $L$-ordered set, where for $\mathcal{A}\in L^{(L^X)}$, $$\sqcup\mathcal{A}=\bigvee\limits_{A\in L^X}A\wedge\mathcal{A}(A),\ \sqcap \mathcal{A}=\bigwedge\limits_{A\in L^X}\mathcal{A}(A)\rightarrow A.$$ If the background set is clear, then we always drop the subscript in ${\rm sub}_X$ to be ${\rm sub}$.
\end{Exam}

In a complete $L$-ordered set $(X,e)$, there are a top element $\top$ and a bottom element $\bot$
such that $e(x,\top)=e(\bot,x)=1$ for every $x\in X$ \cite{YaoFrm}.

\begin{Rem}{\rm(1)} $x=\sqcup A$ if{}f $e(x,y)={\rm sub}({\downarrow}y,A)\ (\forall y\in X)$ if{}f

\ \ \ {\rm(J1)} $A(y)\leq e(y,x)\ (\forall y\in X)$;

\ \ \ {\rm(J2)} $\bigwedge\limits_{z\in X}A(z)\rightarrow e(z,y)\leq e(x,y)\ (\forall y\in X)$.

{\rm(2)} $x=\sqcap A$ if{}f $e(y,x)={\rm sub}(A,{\uparrow}y)\ (\forall y\in X)$  if{}f

\ \ \ {\rm(M1)} $A(y)\leq e(x,y)\ (\forall y\in X)$;

\ \ \ {\rm(M2)} $\bigwedge\limits_{z\in X}A(z)\rightarrow e(y,z)\leq e(y,x)\ (\forall y\in X)$.

{\rm(3)} $A(x)\leq e(x,\sqcup A)\wedge e(\sqcap A,x)\ (\forall x\in X)$, that is, $A\leq{\uparrow}(\sqcap A)\wedge {\downarrow}(\sqcup A)$.

{\rm(4)} Let $f:X\longrightarrow Y$ be a mapping between sets. Then ${\rm sub}(f_L^\rightarrow(A),B)={\rm sub}(A,f_L^\leftarrow(B))$ for all $A\in L^X,\ B\in L^Y$.\end{Rem}

\subsection{Open filters of $L$-valued topological spaces}

\begin{Def}Let $(X,\mathcal{O}(X))$ be an $L$-valued topological space.
An open filter of $X$ is a mapping $u:\mathcal{O}(X)\longrightarrow L$
satisfying that

{\rm(F1)} $u(A\wedge B)=u(A)\wedge u(B)$ $(\forall A,B\in\mathcal{O}(X))$;

{\rm(F2)} $u(a_X)\geq a$ $(\forall a\in L)$.

Let $\Phi_L(X)$ denote the set of all open filters on $X$.\end{Def}

\begin{Exam}Let $(X,\mathcal{O}(X)$ be an $L$-valued topological space.

{\rm(1)} For $A\in L^X$, define $[A]:\mathcal{O}(X)\longrightarrow L$ by $[A](B)={\rm sub}(A,B)\ (\forall B\in\mathcal{O}(X))$. Then $[A]\in\Phi_L(X)$.

{\rm(2)} For $x\in X$, Define $[x]:\mathcal{O}(X)\longrightarrow L$ by $[x](B)=B(x)\ (\forall B\in\mathcal{O}(X))$. Then $[x]\in\Phi_L(X)$, called the {\it pointed filter} of $x$, which is an $L$-valued case of the open neighborhood system of $x$.\end{Exam}

\begin{Pro}Let $X$ be an $L$-valued topological space. For every $u\in\Phi_L(X)$, it holds that: for every $B\in\mathcal{O}(X)$,

{\rm(F3)} $u(B)=\bigvee\limits_{A\in\mathcal{O}(X)}u(A)\wedge {\rm sub}(A,B)$;

{\rm(F3$'$)} $u=\bigvee\limits_{A\in\mathcal{O}(X)}u(A)\wedge [A]$.

{\rm(F4)} $u(B)=\bigwedge\limits_{A\in\mathcal{O}(X)}{\rm sub}(B,A)\rightarrow u(A)$;

{\rm(F4$'$)} $u(B)={\rm sub}([B],u)$.
\end{Pro}

\p{(F3)$\Longleftrightarrow$(F3$'$) is clear. Firstly, $$\bigvee\limits_{A\in\mathcal{O}(X)}u(A)\wedge {\rm sub}(A,B)\geq u(B)\wedge {\rm sub}(B,B)=u(B).$$
Secondly, for every $A\in\mathcal{O}(X)$, it is easily seen that $A\wedge ({\rm sub}(A,B))_X\leq B$ and then $$u(A)\wedge {\rm sub}(A,B)\leq u(A)\wedge u(({\rm sub}(A,B))_X)=u(A\wedge ({\rm sub}(A,B))_X)\leq u(B).$$ Hence, $u(B)=\bigvee\limits_{A\in\mathcal{O}(X)}u(A)\wedge {\rm sub}(A,B)$.

(F4)$\Longleftrightarrow$(F4$'$) is clear. Firstly, $$\bigwedge\limits_{A\in\mathcal{O}(X)}{\rm sub}(B,A)\rightarrow u(A)\leq {\rm sub}(B,B)\rightarrow u(B)=u(B).$$
Secondly, for every $A\in\mathcal{O}(X)$, by the proof of (F3), $u(B)\wedge {\rm sub}(B,A)\leq u(A)$ and then $u(B)\leq {\rm sub}(B,A) \rightarrow u(A)$. Hence, $u(B)=\bigwedge\limits_{A\in\mathcal{O}(X)}{\rm sub}(B,A)\rightarrow u(A)$.}

\begin{Pro}Let $A,B\in L^X$. If $A\leq B$, then ${\rm sub}(B,C)\leq {\rm sub}(A,C)$ for each $C\in L^X$ and so $[B]\leq [A]$.\end{Pro}

\p{This is trivial since ${\rm sub}(A,B)=1$.}

For an $L$-valued topological space $(X,\mathcal{O}(X))$, we can define the convergence of open filters as a mapping $\lim_X:\Phi_L(X)\times X\longrightarrow L$ by $\lim_Xu(x)={\rm sub}_X([x],u)$, that is, $\lim_Xu(x)=\bigwedge\limits_{A\in\mathcal{O}(X)}A(x)\rightarrow u(A)$.

\section{The open filter monad of $L$-valued topological spaces}

In this section will establish the open filter monad for $L$-valued topological spaces.

Let $X$ be an $L$-valued topological space. For every $A\in\mathcal{O}(X)$, define $\phi(A):\Phi_L(X)\longrightarrow L$ by $\phi(A)(u)=u(A)$.
Equipping $\Phi_L(X)$ with an $L$-valued topology generated by $\{\phi(A)|\ A\in\mathcal{O}(X)\}$ as a base, denoted by  $\mathcal{O}(\Phi_L(X))$.
Clearly, $g:\Phi_L(X)\longrightarrow \Phi_L(Y)$ is continuous if{}f $g_L^\leftarrow(\phi(B))\in\mathcal{O}(\Phi_L(X))\ (\forall B\in\mathcal{O}(Y))$.

For $\Phi_L(X)$, $\alpha\in\Phi_L^2(X)$ means an open filter $\alpha:\mathcal{O}(\Phi_L(X))\longrightarrow L$. For $A\in\mathcal{O}(X)$, the mapping $\phi(\phi(A)):\Phi^2_L(X)\longrightarrow L$ is given by $\phi(\phi(A))(\alpha)=\alpha(\phi(A))$.

\begin{Pro}Let $\mathcal{B}$ be a base of an $L$-topological space $(X,\mathcal{O}(X))$. Then For every $A\in\mathcal{O}(X))$ and every $x\in X$, it holds that $$A(x)=\bigvee\limits_{B\in\mathcal{B}}B(x)\wedge {\rm sub}(B,A).$$\end{Pro}

\p{Firstly, for every $B\in\mathcal{B}$ and every $x\in X$, we have $B(x)\wedge {\rm sub}(B,A)\leq A(x)$ and then $A(x)\geq\bigvee\limits_{B\in\mathcal{B}}B(x)\wedge {\rm sub}(B,A).$

Secondly, there exists $\{(B_j,a_j)\}_{j\in J}\subseteq\mathcal{B}\times L$ such that $A=\bigvee_j B_j\wedge(a_j)_X$.
For each $j\in J$, $a_j\leq {\rm sub}(B_j,A)$ and then $(B_j\wedge(a_j)_X)(x)=B_j(x)\wedge a_j\leq B_j(x)\wedge {\rm sub}(B_j,A)$. Thus $A(x)\leq\bigvee\limits_{B\in\mathcal{B}}B(x)\wedge {\rm sub}(B,A).$

Hence, $A(x)=\bigvee\limits_{B\in\mathcal{B}}B(x)\wedge {\rm sub}(B,A).$}

\vskip 2mm

Let $f:X\longrightarrow Y$ be a mapping. Define $\Phi_L f:\Phi_L(X)\longrightarrow \Phi_L(Y)$ by $$(\Phi_L f)(u)(B)=u(f_L^\leftarrow(B)).$$It is easily shown that $\Phi_Lf$ is a well-defined mapping.

\begin{Lem}{\rm(1)} The space $\Phi_L(X)$ is $T_0$.

{\rm(2)} If $f:X\longrightarrow Y$ is a continuous mapping, then so is $\Phi_L f:\Phi_L(X)\longrightarrow \Phi_L Y$.

Consequently, $\Phi_L$  can be considered as an endofunctor on $L$-${\bf Top}_0$.\end{Lem}

\p{(1) Let $u,v\in\Phi_L(X)$ and $W(u)=W(v)$ for all $W\in\mathcal{O}(\Phi_L(X))$. Then $\phi(A)(u)=\phi(A)(v)$ for all $A\in\mathcal{O}(X)$. That is, $u(A)=v(A)$ for all $A\in\mathcal{O}(X)$, and so $u=v$. Therefore, $\Phi_L(X)$ is $T_0$.

(2) For every $B\in\mathcal{O}(Y)$ and every $u\in\Phi_L(X)$. Then
$$\begin{array}{lll}(\Phi_L f)_L^\leftarrow(\phi(B))(u)&=&(\phi(B))((\Phi_L f)(u))\\ \ &=&(\Phi_L f)(u)(B)\\ \ &=&u(f_L^\leftarrow(B))\\ \ &=&\phi(f_L^\leftarrow(B))(u).\end{array}$$Therefore, $(\Phi_L f)_L^\leftarrow(\phi(B))=\phi(f_L^\leftarrow(B))\in\mathcal{O}(\Phi_L(X))$. Hence, $\Phi_L f:\Phi_L(X)\longrightarrow \Phi_L(Y)$ is continuous.}

\vskip 2mm

Define $\eta_X:X\longrightarrow\Phi_L(X)$ by $\eta_X(x)=[x]\ (\forall x\in X)$.

\begin{Lem}Let $X$ be an $L$-valued topological space. Then for each $U\in\mathcal{O}(X)$, it holds that $(\eta_X)_L^\leftarrow(\phi(U))=U$.\end{Lem}

\p{For each $x\in X$, $$(\eta_X)_L^\leftarrow(\phi(U))(x)=\phi(U)(\eta_X(x))=\phi(U)([x])=[x](U)=U(x),$$as desired.}

\begin{Th}$\eta:id_{L\mbox{-}{\bf Top}_0}\longrightarrow \Phi_L$ is a natural transformation.\end{Th}

\p{(1) For every $U\in\mathcal{O}(X)$, by Lemma 3.3, we have $(\eta_X)_L^\leftarrow(\phi(U))=U\in\mathcal{O}(X)$ and then $\eta_X:X\longrightarrow\Phi_L(X)$ is continuous.

(2) Let $f:X\longrightarrow Y$ be a continuous mapping. We need to show that $\eta_Y\cdot f=(\Phi_Lf)\cdot \eta_X$. In fact, for every $x\in X$ and every $B\in\mathcal{O}(Y)$, we have
$(\Phi_L f)\cdot \eta_X(x)(B)=(\Phi_L f)([x])(B)=[x](f_L^\leftarrow(B))=f_L^\leftarrow(B)(x)=B(f(x))=[f(x)](B)=\eta_Y\cdot f(x)(B).$}

\vskip 2mm

Let $\alpha\in\Phi^2_L(X)$, define $\mu_X(\alpha)(A)=\alpha(\phi(A))\ (\forall A\in\mathcal{O}(X))$.

\begin{Th}$\mu_X(\alpha)\in\Phi_L(X)$ for every $\alpha\in\Phi^2_L(X)$.\end{Th}

\p{For all $A,B\in\mathcal{O}(X)$, $$\begin{array}{lll}\mu_X(\alpha)(A)\wedge\mu_X(\alpha)(B)&=&\alpha(\phi(A))\wedge\alpha(\phi(B))\\\ &=&\alpha(\phi(A)\wedge\phi(B))\\ \ &=&\alpha(\phi(A\wedge B))\\ \ &=&\mu_X(\alpha)(A\wedge B).\end{array}$$ For all $a\in L$, $\mu_X(\alpha)(a_X)=\alpha(\phi(a_X))\geq \alpha(a_{\Phi_L(X)})\geq a$. Notice here $\phi(a_X)(u)=u(a_X)\geq a\ (\forall u\in\Phi_L(X))$ and then $\phi(a_X)\geq a_{\Phi_L(X)}$.}

\begin{Th}$\mu:\Phi^2_L\longrightarrow\Phi_L$ is a natural transformation.\end{Th}

\p{Step 1. For every $L$-valued topological space $X$ and every $A\in\mathcal{O}(X)$, we have  $$\begin{array}{lll}(\mu_X)^\leftarrow_L(\phi(A))(\alpha)&=&\phi(A)(\mu_X(\alpha))\\ \ &=&\mu_X(\alpha)(A)\\ \ &=&\alpha(\phi(A))\\ \ &=&\phi\phi(A)(\alpha)\\ \ &\in&\mathcal{O}(\Phi^2_L(X)).\end{array}$$Hence, $\mu_X:\Phi^2_L(X)\longrightarrow\Phi_L(X)$ is a continuous mapping.

Step 2. Let $f:X\longrightarrow Y$ be a mapping. We need to prove that $\mu_Y\cdot(\Phi^2_L f)=\Phi_L f\cdot \mu_X$. For every $\alpha\in \Phi^2_L(X)$, for every $B\in\mathcal{O}(Y)$,
for every $u\in\Phi_L(Y)$, by the proof of Lemma 3.2(2), we have $(\Phi_Lf)_L^\leftarrow(\phi(B))=\phi(f_L^\leftarrow(B))$, and then
$$\begin{array}{lll}(\mu_Y)\cdot(\Phi^2_L f)(\alpha)(B)&=&(\Phi_L(\Phi_Lf))(\alpha)(\phi(B))\\ \ &=&\alpha((\Phi_L f)_L^\leftarrow(\phi(B)))\\ \ &=&\alpha(\phi(f_L^\leftarrow(B)))\\ \ &=&\mu_X(\alpha)(f_L^\leftarrow(B))\\ \ &=&\Phi_L f\cdot \mu_X(\alpha)(B).\end{array}$$ Hence, $\mu_Y\cdot(\Phi^2_L f)=\Phi_L f\cdot \mu_X$.}

\begin{Th}$(\Phi_L,\eta,\mu)$ is a monad over $L\mbox{-}{\bf Top}_0$.\end{Th}

\p{(1) $\mu\circ (\Phi_L\mu)=\mu\circ (\mu\Phi_L)$, that is, $\mu\cdot (\Phi_L\mu_X)=\mu\cdot (\mu_{\Phi_L(X)})$ for every $T_0$ $L$-valued topological space, where $\Phi_L\mu_X:\Phi_L(\Phi^2_L(X))\longrightarrow \Phi_L(\Phi_L(X))$ and $\mu_{\Phi_L(X)}:\Phi^2_L(\Phi_L(X))\longrightarrow\Phi_L(\Phi_L(X))$.
Notice that for every $\alpha\in\Phi^2_L(X)$, $$(\mu_X)_L^\leftarrow(\phi(U))(\alpha)=\phi(U)(\mu_X(\alpha))=\mu_X(\alpha)(U)=\phi\phi(U)(\alpha).$$ Then  $(\mu_X)_L^\leftarrow(\phi(U))=\phi\phi(U)$ and for every $\Xi\in\Phi_L(\Phi^2_L(X))$, we have
$$\begin{array}{lll}\mu[(\Phi_L\mu_X)(\Xi)](U)&=&\mu[(\mu_X)_L^\rightarrow(\Xi)](U)\\ \ &=&(\mu_X)_L^\rightarrow(\Xi)(\phi(U))
\\ \ &=&\Xi((\mu_X)_L^\leftarrow(\phi(U)))\\ \ &=&\Xi(\phi\phi(U))\\ \ &=&(\mu_{\Phi X})(\Xi)(\phi(U))\\ \ &=&\mu[(\mu_{\Phi X})(\Xi)](U).\end{array}$$Hence, $\mu\cdot (\Phi_L\mu_X)=\mu\cdot (\mu_{\Phi_L(X)})$.

(2) $\mu\circ\eta\Phi_L=id_{\Phi_L}=\mu\circ \Phi_L\eta$, that is, $\mu(\eta_{\Phi_L(X)})(u)=u=\mu(\Phi_L\eta)(u)$ for each $L$-valued topological space $X$ and every $u\in \Phi_L(X)$. In fact, for all $u\in\Phi_L(X)$ and $U\in\mathcal{O}(X)$, we have
$$\mu(\eta_{\Phi_L(X)}(u))(U)=\mu([u])(U)=[u](\phi(U))=\phi(U)(u)=u(U)$$and by Lemma 3.3,
$$\mu[\Phi_L\eta(u)](U)=\mu[(\eta_X)_L^\rightarrow(u)](U)=(\eta_X)_L^\rightarrow(u)(\phi(U))
=u[(\eta_X)_L^\leftarrow(\phi(U))]=u(U).$$ Hence, $\mu\circ\eta\Phi_L=id_{\Phi_L}=\mu\circ \Phi_L\eta$.}

\section{$L$-valued continuous lattices and $L$-valued Scott topology}

In order to characterize $L$-valued domain structures via the monad $(\Phi_L,\eta,\mu)$,
we will recall some basic concepts and results about $L$-valued domains and $L$-valued Scott topology \cite{PartI,PartII}.

\begin{Def}Let $X$ be an $L$-ordered set. An $L$-subset $D\in L^X$ is called directed if

{\rm(D1)} $\bigvee\limits_{x\in X}D(x)=1$;

{\rm(D2)} $D(x)\wedge D(y)\leq\bigvee\limits_{z\in X}D(z)\wedge e(x,z)\wedge e(y,z)$.

An $L$-subset $I\in L^X$ is called an {\it ideal} of $X$ if it is a directed lower set. Denote by $\mathcal{D}_L(X)$ (resp., $\ mathcal{I}_L(X)$) the set of all directed $L$-subsets (resp., ideals) of $X$. An $L$-ordered set $X$ is called an $L$-{\it valued dcpo} if every directed $L$-subset has a supremum, or equivalently, every ideal has a supremum.
\end{Def}

\begin{Def}For an $L$-valued dcpo $X$, $A\in L^X$ is called {\it Scott open} if it satisfies one of the following equivalent conditions:

{\rm(1)} $A(\sqcup D)=\sqcup A_L^\rightarrow(D)$ for every $D\in\mathcal{D}_L(X)$;

{\rm(2)} $A(\sqcup I)=\sqcup A_L^\rightarrow(I)$ for every $I\in\mathcal{I}_L(X)$;

{\rm(3)} $A$ is an upper set and $A(\sqcup D)\leq\sqcup A_L^\rightarrow(D)$ for every $D\in\mathcal{D}_L(X)$;

{\rm(4)} $A$ is an upper set and $A(\sqcup I)\leq\sqcup A_L^\rightarrow(I)$ for every $I\in\mathcal{I}_L(X)$.

The family $\sigma_L(X)$ of all Scott open $L$-subsets of $X$ forms an $L$-valued topology, called the $L$-{\it valued Scott topology} on $X$.\end{Def}

\begin{Def}Let $(X,e)$ be a complete $L$-ordered set. For every $x\in X$, define ${\Downarrow}x\in L^X$ by $${\Downarrow}x(y)=\bigwedge\limits_{I\in\mathcal{I}_L(X)}e(x,\sqcup I)\rightarrow I(y).$$ A complete $L$-ordered set $(X,e)$ is called an $L$-valued continuous lattice
if $\sqcup{\Downarrow}x=x$ for each $x\in X$.\end{Def}

The $L$-relation ${\Downarrow}$ is an $L$-valued version of the classical way-below relation.

\begin{Pro}
{\rm(1)} ${\Downarrow}x(y)\leq e(y,x)\ (\forall x,y\in X)$;

{\rm(2)} $e(y_1,y)\wedge{\Downarrow}x(y)\wedge e(x,x_1)\leq{\Downarrow}w(z)\ (\forall x,x_1,y,y_1\in X)$.
\end{Pro}

\begin{Pro}If $(X,e)$ is an $L$-valued continuous lattice, then

{\rm(1)} The mapping ${\Downarrow}$ is interpolative, that is, ${\Downarrow}x(y)=\bigvee\limits_{z\in X}{\Downarrow}z(y)\wedge {\Downarrow}x(z)$;

{\rm(2)} $\{{\Uparrow}x|\ x\in X\}$ is a base of $\sigma_L(X)$, where ${\Uparrow}x(y)={\Downarrow}y(x)\ (\forall x,y\in X)$;

{\rm(3)} $\sigma_L(X)$ is a $T_0$ $L$-valued topological space.\end{Pro}

\begin{Lem}\label{wbI} Let $I\in \mathcal{I}_L(X)$ and $x\in X$. Then

{\rm(1)} If $x=\sqcup I$, then ${\Downarrow}x\leq I$;

{\rm(2)} ${\Downarrow}(\sqcup I)(x)\leq I(x)$.\end{Lem}

\p{(1) For each $y\in X$, ${\Downarrow}x(y)=\bigwedge\limits_{J\in\mathcal{I}_L(X)}e(x,\sqcup J)\rightarrow J(y)\leq e(x,\sqcup I)\rightarrow I(y)=1\rightarrow I(y)=I(y).$

(2) ${\Downarrow}(\sqcup I)(x)=\bigwedge\limits_{J\in\mathcal{I}_L(X)}e(\sqcup I,\sqcup J)\rightarrow J(x)\leq e(\sqcup I,\sqcup I)\rightarrow I(x)=I(x).$}

Let $X$ be a complete $L$-ordered set and equipped with the $L$-valued Scott topology $\sigma_L(X)$. For every $u\in\Phi_L(X)$, define
$u^l\in L^X$ by $$u^l(x)=\bigvee\limits_{A\in\mathcal{O}(X)}u(A)\wedge A^l(x),$$ and further define
$\lim_S:\Phi_L(X)\times X\longrightarrow L$ by $$\lim{}_S u(x)=\bigvee\limits_{I\in\mathcal{I}_L(X)}{\rm sub}_X(I,u^l)\wedge e(x,\sqcup I).$$ This is the fuzzy Scott convergence of open filters, which is a modification of Scott convergence of stratified $L$-filters in [Yao].

\begin{Lem}Let $X$ be a complete $L$-ordered set and $u$ be an open filter of $\sigma_L(X)$.  Then $u^l\in\mathcal{I}_L(X)$.\end{Lem}

\p{Firstly, it is routine to show that $u^l$ is a lower set.

Secondly, $\bigvee\limits_{x\in X}u^l(x)\geq u^l(\bot)=\bigvee\limits_{A\in\mathcal{O}(X)}u(A)\wedge {\rm sub}(A,{\uparrow}\bot)\geq u(1_X)=1$.

Thirdly, for all $x,y\in X$, $$\begin{array}{lll}u^l(x)\wedge u^l(y)
&=&\bigvee\limits_{A\in\sigma_L(X)}u(A)\wedge {\rm sub}(A,{\uparrow}x)\wedge\bigvee\limits_{B\in\sigma_L(X)}u(B)\wedge {\rm sub}(B,{\uparrow}y)
\\ &=&\bigvee\limits_{A,B\in\sigma_L(X)}u(A)\wedge u(B)\wedge {\rm sub}(A,{\uparrow}x)\wedge {\rm sub}(B,{\uparrow}y)
\\&\leq&\bigvee\limits_{A,B\in\sigma_L(X)}u(A\wedge B)\wedge {\rm sub}(A\wedge B,{\uparrow}x)\wedge {\rm sub}(A\wedge B,{\uparrow}y)
\\&\leq& \bigvee\limits_{C\in\sigma_L(X)}u(C)\wedge {\rm sub}(C,{\uparrow}x)\wedge {\rm sub}(C,{\uparrow}y)
\\&\leq&\bigvee\limits_{z\in X}\bigvee\limits_{C\in\sigma_L(X)}u(C)\wedge {\rm sub}(C,{\uparrow}z)\wedge e(x,z)\wedge e(y,z).\end{array}$$

Hence, $u^l\in\mathcal{I}_L(X)$}

For every open filter $u$ of $X$ with respect to the $L$-valued Scott topology. Define $r:\Phi_L(X)\longrightarrow X$ by $r(u)=\sqcup u^l$.

\begin{Th}$\lim_S u(x)=e(x,r(u))$.\end{Th}

\p{On one hand, $$\lim{}_S u(x)=\bigvee\limits_{I\in\mathcal{I}_L(X)}{\rm sub}_X(I,u^l)\wedge e(x,\sqcup I)\geq {\rm sub}_X(u^l,u^l)\wedge e(x,\sqcup u^l)=e(x,\sqcup u^l).$$On the other hand $$\lim{}_S u(x)=\bigvee\limits_{I\in\mathcal{I}_L(X)}{\rm sub}_X(I,u^l)\wedge e(x,\sqcup I)\leq \bigvee\limits_{I\in\mathcal{I}_L(X)}e(\sqcup I,\sqcup u^l)\wedge e(x,\sqcup I)\leq e(x,\sqcup u^l).$$
Hence, $\lim_S u(x)=e(x,r(u))$.}

\begin{Rem}The meaning of Theorem 4.8 is that, the set of limit point of an open filter $u$ is a lower set with $r(u)$ as the largest element.\end{Rem}

\begin{Pro}Let $X$ be a complete $L$-ordered set. Then $X$ is an $L$-valued continuous lattice if{}f $x=\sqcup[x]^l$ holds for every $x\in X$.\end{Pro}

\p{$\Longrightarrow$: On one hand, for every $y\in X$, $$[x]^l(y)=\bigvee\limits_{A\in\mathcal{O}(X)}[x](A)\wedge {\rm sub}(A,{\uparrow}y)\leq A(x)\wedge(A(x)\rightarrow e(y,x))\leq e(y,x).$$ On the other hand,
$$\begin{array}{lll}\ &\bigwedge\limits_{z\in X}[x]^l(z)\rightarrow e(x,y)
\\=&\bigwedge\limits_{z\in X}\bigwedge\limits_{A\in\sigma_L(X)}([x](A)\wedge {\rm sub}(A,{\uparrow}z))\rightarrow e(z,y)&({\rm Prop. 2.2(6)})
\\\leq& \bigwedge\limits_{z\in X}\bigwedge\limits_{x_1\in X}({\Uparrow}x_1(x)\wedge {\rm sub}({\Uparrow}x_1,{\uparrow}z))\rightarrow e(z,y)&({\rm Prop. 4.5(2)})
\\=&\bigwedge\limits_{x_1\in X}{\Downarrow}x(x_1)\rightarrow[\bigwedge\limits_{z\in X}{\rm sub}({\Uparrow}x_1,{\uparrow}z)\rightarrow e(z,y)]&({\rm Prop. 2.2(7,12)})
\\\leq& \bigwedge\limits_{x_1\in X}{\Downarrow}x(x_1)\rightarrow[\bigwedge\limits_{z\in X}{\rm sub}({\uparrow}x_1,{\uparrow}z)\rightarrow e(z,y)]&({\rm Prop. 4.4(1)})
\\\leq& \bigwedge\limits_{x_1\in X}{\Downarrow}x(x_1)\rightarrow[\bigwedge\limits_{z\in X}e(z,x_1)\rightarrow e(z,y)]
\\\leq& \bigwedge\limits_{x_1\in X}{\Downarrow}x(x_1)\rightarrow e(x_1,y)&({\rm Prop. 2.7(3)})
\\=&e(\sqcup{\Downarrow}x,y)=e(x,y).\end{array}$$

Hence, $x=\sqcup[x]^l$.

$\Longleftarrow$: By Proposition 2.7(1), we only need to show that $[x]^l\leq{\Downarrow}x\leq {\downarrow}x$. Firstly, ${\Downarrow}x\leq {\downarrow}x$ is obvious by Proposition 4.4(1). Secondly,
since ${\Downarrow}x(y)=\bigwedge\limits_{I\in\mathcal{I}_L(X)}e(x,\sqcup I)\rightarrow I(y)$ and $[x]^l(y)=\bigvee\limits_{A\in\sigma_L(X)}A(x)\wedge {\rm sub}(A,{\uparrow}y)$, we have
$$\begin{array}{lll}\ &A(x)\wedge {\rm sub}(A,{\uparrow}y)\wedge e(x,\sqcup I)
\\\leq& A(\sqcup I)\wedge{\rm sub}(A,{\uparrow}y)&({\rm Def.4.2})
\\=&\bigvee\limits_{z\in X}A(z)\wedge I(z)\wedge {\rm sub}(A,{\uparrow}y)&({\rm Def.4.2,\ Rem.2.9(2)})
\\\leq&\bigvee\limits_{z\in X} I(z)\wedge {\uparrow}y(z)
\\ =& \bigvee\limits_{z\in X} I(z)\wedge e(y,z)
\\\leq& I(y).&({\rm Def.4.1})\end{array}$$Then $A(x)\wedge {\rm sub}(A,{\uparrow}y)\leq e(x,\sqcup I)\rightarrow I(y)$ and consequently, $[x]^l\leq{\Downarrow}x$ by the arbitrariness of $y\in X$.}

\section{The first main theorem}

{\bf The First Main Theorem}\ \ \ \ If $X$ is an $L$-valued continuous lattice and equipping $X$ with the $L$-valued Scott topology,
then the pair $(X,r)$ is a $\Phi_L$-algebra over the monad $(\Phi_L,\mu,\eta)$, that is, $r\cdot\Phi_L r=r\cdot \mu_X$ and $r\cdot\eta_X=id_X$, where $r:\Phi_L(X)\longrightarrow L$ is given by $r(u)=\sqcup u^l$.

First of all, by Proposition 4.5(3), $\sigma_L(X)$ is a $T_0$ $L$-valued topological space.

\begin{Lem}Let $(X,e)$ be a complete $L$-ordered set. For every $A\in \sigma_L(X)$, it holds that $r^\leftarrow_L(A)\leq\phi(A)$.\end{Lem}

\p{For all $u\in\Phi_L(X)$, we have
$$\begin{array}{llll}r^\leftarrow_L(A)(u)
&=&A(r(u))=A(\sqcup u^l)=\sqcup A_L^\rightarrow(u^l)&({\rm Def. 4.2})
\\\ &=&\bigvee\limits_{x\in X}A(x)\wedge u^l(x)&({\rm Exam. 2.9(2)})
\\ \ &=&\bigvee\limits_{x\in X}A(x)\wedge\bigvee\limits_{B\in\sigma_L(X)}u(B)\wedge{\rm sub}(B,{\uparrow}x)
\\ \ &=&\bigvee\limits_{x\in X}\bigvee\limits_{B\in\sigma_L(X)}{\rm sub}({\uparrow}x,A)\wedge u(B)\wedge {\rm sub}(B,{\uparrow}x)
\\ \ &\leq&\bigvee\limits_{B\in\sigma_L(X)} u(B)\wedge {\rm sub}(B,A)&({\rm Exam. 2.9(3)})
\\ \ &=& u(A)=\phi(A)(u).&({\rm F3})\end{array}$$
Hence, $r^\leftarrow_L(A)\leq \phi(A)$. }

\begin{Lem}Let $(X,e)$ be an $L$-valued continuous lattice. For every $A\in\sigma_L(X)$ and every $x\in X$, it holds that $${\Downarrow}(\sqcap A)(x)\leq{\rm sub}(\phi(A),r_L^\leftarrow({\Uparrow}x)).$$\end{Lem}

\p{For every $u\in\Phi_L(X)$, we have
$$\begin{array}{llll}e(\sqcap A,r(u))&=&e(\sqcap A,\sqcup u^l)\\ \ &\geq& u^l(\sqcap A)&({\rm M1})
\\ \ &=&\bigvee\limits_{B\in \sigma_L(X)}u(B)\wedge{\rm sub}(B,{\uparrow}(\sqcap A))
\\ \ &\geq& u(A)\wedge{\rm sub}(A,{\uparrow}(\sqcap A))
\\ \ &=&u(A)=\phi(A)(u).&({\rm Rem. 2.10(3)})\end{array}$$
Then for each $x\in X$, by Proposition 4.4(2), we have
$$\begin{array}{lll}\phi(A)(u)\wedge {\Downarrow}(\sqcap A)(x)&\leq& e(\sqcap A,r(u))\wedge  {\Downarrow}(\sqcap A)(x)\\ \ &\leq& {\Downarrow}(r(u))(x)={\Uparrow}x(r(u))=r_L^\leftarrow({\Uparrow}x)(u).\end{array}$$
Hence,
${\Downarrow}(\sqcap A)(x)\leq \bigwedge\limits_{u\in\Phi_L(X)}\phi(A)(u)\rightarrow r_L^\leftarrow({\Uparrow}x)(u)={\rm sub}(\phi(A),r_L^\leftarrow({\Uparrow}x))$.}

\begin{Pro}$r\cdot\Phi_L r=r\cdot \mu_X$.\end{Pro}

\p{Let $\alpha\in\Phi^2_L(X)$. Then
$$r[\Phi_L r(\alpha)]=\sqcup (\Phi_L r(\alpha))^l,\ \ \ r[\mu_X(\alpha)]=\sqcup (\mu_X(\alpha))^l.$$
For every $x\in X$,
$$\begin{array}{lll}(\Phi_L r(\alpha))^l(x)
&=&\bigvee\limits_{A\in \sigma_L(X)}\Phi_L r(\alpha)(A)\wedge {\rm sub}(A,{\uparrow}x)
\\\ &=&\bigvee\limits_{A\in \sigma_L(X)}r_L^\rightarrow(\alpha)(A)\wedge {\rm sub}(A,{\uparrow}x)
\\ \ &=&\bigvee\limits_{A\in\sigma_L(X)} \alpha(r_L^\leftarrow(A))\wedge {\rm sub}(A,{\uparrow}x)\end{array}$$
and
$$\begin{array}{lll}(\mu_X(\alpha))^l(x)
&=&\bigvee\limits_{A\in \sigma_L(X)}\mu_X(\alpha)(A)\wedge {\rm sub}(A,{\uparrow}x)
\\ \ &=&\bigvee\limits_{A\in\sigma_L(X)}\alpha(\phi(A))\wedge{\rm sub}(A,{\uparrow}x).\end{array}$$

Firstly, by Lemma 5.1, $(\Phi_L r(\alpha))^l\leq(\mu_X(\alpha))^l$. Hence, $$r[\mu_X(\alpha)]=\sqcup(\Phi_L r(\alpha))^l\leq\sqcup(\mu_X(\alpha))^l=r[\Phi r(\alpha)].$$

Secondly, let $x_0=r[\mu_X(\alpha)]$. Then $x_0=\sqcup{\Downarrow}x_0$ and by Lemma \ref{wbI}(1), ${\Downarrow}x_0\leq (\mu_X(\alpha))^l$. We will show that
${\Downarrow}x_0\leq (\Phi_L r(\alpha))^l$, so that $r[\mu_X(\alpha)]=x_0\leq \sqcup (\Phi_L r(\alpha))^l=r[\Phi r(\alpha)]$.

Then
$$\begin{array}{llll}{\Downarrow}x_0(x)
&=&\bigvee\limits_{y\in X}{\Downarrow}x_0(y)\wedge{\Downarrow}y(x)&({\rm Prop. 4.5(1)})
\\ \ &\leq& \bigvee\limits_{y\in X}(\mu_X(\alpha))^l(y)\wedge{\Downarrow}y(x)
\\ \ &=&\bigvee\limits_{y\in X}\bigvee\limits_{A\in \sigma_L(X)}\alpha(\phi(A))\wedge {\rm sub}(A,{\uparrow}y)\wedge{\Downarrow}y(x)
\\ \ &=& \bigvee\limits_{y\in X}\bigvee\limits_{A\in \sigma_L(X)}\alpha(\phi(A))\wedge e(y,\sqcap A)\wedge{\Downarrow}y(x)&({\rm Rem. 2.10(2)})
\\ \ &\leq& \bigvee\limits_{A\in\sigma_L(X)}\alpha(\phi(A))\wedge{\Downarrow}(\sqcap A)(x)&({\rm Prop. 4.4(2)})
\\ \ &\leq& \bigvee\limits_{A\in \sigma_L(X)}\alpha(\phi(A))\wedge {\rm sub}(\phi(A),r_L^\leftarrow({\Uparrow}x))&({\rm Lem. 5.2})
\\ \ &\leq& \alpha(r_L^\leftarrow({\Uparrow}x))&({\rm F4})
\\ \ &\leq& \bigvee\limits_{A\in \sigma_L(X)}\alpha(r_L^\leftarrow(A))\wedge{\rm sub}(A,{\uparrow}x)&({\rm Prop. 4.4(1)})
\\ \ &=&  (\Phi_L r(\alpha))^l(x).\end{array}$$
Then ${\Downarrow}x_0\leq(\Phi_L r(\alpha))^l$ and
$$r[\Phi_L r(\alpha)]=x_0=\sqcup{\Downarrow}x_0\leq\sqcup(\Phi_L r(\alpha))^l=r[\mu_X(\alpha)].$$

Therefore, $r\cdot\Phi_L r=r\cdot \mu_X$.}

\begin{Pro}$r\cdot\eta_X=id_X$.\end{Pro}

\p{For each $x\in X$, we have $$[x]^l(y)=\bigvee\limits_{A\in\sigma_L(X)}[x](A)\wedge {\rm sub}(A,{\uparrow}y)=\bigvee\limits_{A\in\sigma_L(X)}A(x)\wedge{\rm sub}(A,{\uparrow}y).$$
Firstly, $$[x]^l(y)=\bigvee\limits_{A\in\sigma_L(X)}A(x)\wedge {\rm sub}(A,{\uparrow}y)\leq{\uparrow}y(x)={\downarrow}x(y).$$
Secondly, by Proposition 4.5(2), $$[x]^l(y)=\bigvee\limits_{A\in\sigma_L(X)}A(x)\wedge {\rm sub}(A,{\uparrow}y)\geq{\Uparrow}y(x)\wedge {\rm sub}({\Uparrow}y,{\uparrow}y)={\Downarrow}x(y).$$Thus, ${\Downarrow}x\leq[x]^l\leq{\downarrow}x$. Since $\sqcup{\Downarrow}x=x=\sqcup{\downarrow}x$, by Proposition 2.7(1), we have $r\cdot\eta_X(x)=r([x])=\sqcup[x]^l=x$. Hence, $r\cdot\eta_X=id_X$.}

\section{The second main theorem}

{\bf The Second Main Theorem}\ \ \ \  If $(X,r)$ is a $\Phi_L$-algebra over $L$-${\bf Top}_0$,
then by considering $X$ with the specialization $L$-order, $X$ is an $L$-valued continuous lattice and $r(u)=\sqcup u^l$.

\begin{Lem}For every $A\in\mathcal{O}(X)$, it holds that $A\leq r_L^\rightarrow(\phi(A))$.\end{Lem}

\p{For every $x\in X$, we have $$r_L^\rightarrow(\phi(A))(x)=\bigvee\limits_{r(u)=x}\phi(A)(u)=\bigvee\limits_{r(u)=x}u(A)\geq [x](A)=A(x).$$Hence, $A\leq r_L^\rightarrow(\phi(A))$.}

In the following, we will study the property of the specialization $L$-order $e_{\mathcal{O}(X)}$ of the $L$-valued topological space $X$. For simplicity, we write $e$ instead of $e_{\mathcal{O}(X)}$.

\begin{Lem}Let $u,v\in\Phi_L(X)$. Then

{\rm(1)} $\bigwedge\limits_{W\in\mathcal{O}(\Phi_L(X))}W(u)\rightarrow W(v)=\bigwedge\limits_{A\in\mathcal{O}(X)}u(A)\rightarrow v(A)$;

{\rm(2)} ${\rm sub}(u,v)\leq e(r(u),r(v))$.\end{Lem}

\p{(1) Firstly, $$\bigwedge\limits_{W\in\mathcal{O}(\Phi_L(X))}W(u)\rightarrow W(v)\leq \bigwedge\limits_{A\in\mathcal{O}(X)}\phi(A)(u)\rightarrow \phi(A)(v)=\bigwedge\limits_{A\in\mathcal{O}(X)}u(A)\rightarrow v(A).$$Secondly,
$$\begin{array}{llll}\ &\bigwedge\limits_{W\in\mathcal{O}(\Phi_L(X))}W(u)\rightarrow W(v)
\\ \ =&\bigwedge\limits_{\{(A_j,a_j)\}_{j\in J}\subseteq\mathcal{O}(X)\times L}(\bigvee_j\phi(A_j)\wedge (a_j)_X)(u)\rightarrow (\bigvee_j\phi(A_j)\wedge (a_j)_X)
\\ \ \geq&\bigwedge\limits_{\{(A_j,a_j)\}_{j\in J}\subseteq\mathcal{O}(X)\times L}\bigwedge_j(\phi(A_j)\wedge (a_j)_X)(u)\rightarrow (\phi(A_j)\wedge (a_j)_X)(v)&({\rm Prop. 2.2(8)})
\\ \ =&\bigwedge\limits_{\{(A_j,a_j)\}_{j\in J}\subseteq\mathcal{O}(X)\times L}\bigwedge_j(u(A_j)\wedge a_j)\rightarrow (v(A_j)\wedge a_j)
\\ \ \geq&\bigwedge\limits_{\{(A_j,a_j)\}_{j\in J}\subseteq\mathcal{O}(X)\times L}\bigwedge_ju(A_j)\rightarrow v(A_j)&({\rm Prop. 2.2(13)})
\\ \ =&\bigwedge\limits_{A\in\mathcal{O}(X)}u(A)\rightarrow v(A).\end{array}$$
Hence, $\bigwedge\limits_{W\in\mathcal{O}(\Phi_L(X))}W(u)\rightarrow W(v)=\bigwedge\limits_{A\in\mathcal{O}(X)}u(A)\rightarrow v(A)$.

(2) For all $u,v\in\Phi_L(X)$, $$\begin{array}{lll}e(r(u),r(v))&=&\bigwedge\limits_{A\in\mathcal{O}(X)}A(r(u))\rightarrow A(r(v))\\ \ &=&\bigwedge\limits_{A\in\mathcal{O}(X)}r^\leftarrow_L(A)(u)\rightarrow r^\leftarrow_L(A)(v)\\ \ &\geq&\bigwedge\limits_{W\in\mathcal{O}(\Phi_L(X))}W(u)\rightarrow W(v)\\ \ &=&\bigwedge\limits_{A\in\mathcal{O}(X)}u(A)\rightarrow v(A)\\ \ &=&{\rm sub}(u,v).\end{array}$$The proof is completed.}

\begin{Pro}$\lim_Xu(x)=e(x,r(u))$.\end{Pro}

\p{Firstly, by Lemma 6.2(2), $$\lim{}_Xu(x)={\rm sub}([x],u)\leq e(r([x]),r(u))=e(x,r(u)).$$

Secondly, $$e(x,r(u))=\bigwedge\limits_{A\in\mathcal{O}(X)}A(x)\rightarrow A(r(u))=\bigwedge\limits_{A\in\mathcal{O}(X)}A(x)\rightarrow r^\leftarrow_L(A)(u).$$
Since $r_L^\leftarrow(A)\in\mathcal{O}(\Phi_L(X))$ and $u\in\Phi_L(X)$, by applying Proposition 3.1 to the space $(\Phi_L(X),\mathcal{O}(\Phi_L(X)))$
$$\begin{array}{llll}r^\leftarrow_L(A)(u)&=&\bigvee\limits_{B\in\mathcal{O}(X)}\phi(B)(u)\wedge {\rm sub}(\phi(B),r_L^\leftarrow(A))&({\rm Lem. 6.3})
\\ &=&\bigvee\limits_{B\in\mathcal{O}(X)}u(B)\wedge {\rm sub}(r^\rightarrow_L(\phi(B)),A)&({\rm Rem. 2.10(4)})
\\ &\leq&\bigvee\limits_{B\in\mathcal{O}(X)}u(B)\wedge  {\rm sub}(B,A)&({\rm Lem. 6.1})
\\ &=&u(A)&({\rm F3}).\end{array}$$
Hence, $e(x,r(u))\leq \bigwedge\limits_{A\in\mathcal{O}(X)}A(x)\rightarrow u(A)=\lim_Xu(x)$.}

\begin{Pro}The pair $(X,e)$ is a complete $L$-ordered set, where $\sqcap A=r([A])$.\end{Pro}

\p{Let $A\in L^X$. Firstly, by Lemma 6.2(2), $$A(x)\leq{\rm sub}([A],[x])\leq e(r([A]),r([x]))=e(r([A]),x).$$

Secondly, for every $y\in X$,
$$\begin{array}{lll}\ &\bigwedge\limits_{x\in X}A(x)\rightarrow e(y,x)
\\=&\bigwedge\limits_{x\in X}A(x)\rightarrow(\bigwedge\limits_{B\in\mathcal{O}(X)}B(y)\rightarrow B(x))
\\=&\bigwedge\limits_{B\in\mathcal{O}(X)}B(y)\rightarrow(\bigwedge\limits_{x\in X} A(x)\rightarrow B(x))&({\rm Prop. 2.2(7,12)})
\\=&\bigwedge\limits_{B\in\mathcal{O}(X)}B(y)\rightarrow  {\rm sub}(A,B)
\\=&\bigwedge\limits_{B\in\mathcal{O}(X)}[y](B)\rightarrow  [A](B)
\\=& {\rm sub}([y],[A])
\\\leq& e(r([y]),r([A]))&({\rm Lem. 6.2(2)})
\\=&e(y,r([A])).\end{array}$$Thus, $\sqcap A=r([A])$.
Hence, $X$ is complete $L$-ordered set.}

\begin{Lem}The pair $(\Phi_L(X),{\rm sub})$ is an $L$-valued dcpo, where for each directed $L$-subset $\mathcal{A}$  of $\Phi_L(X)$, $\sqcup\mathcal{A}=\bigvee\limits_{u\in\Phi_L(X)}\mathcal{A}(u)\wedge u$.\end{Lem}

\p{Let $\mathcal{A}$ be a directed $L$-subset of $(\Phi_L(X),{\rm sub})$. Put $$w=\bigvee\limits_{u\in\Phi_L(X)}\mathcal{A}(u)\wedge u.$$ We need to show that for every $v\in\Phi_L(X)$, it holds that $w\in\Phi_L(X)$ and $${\rm sub}(w,v)=\bigwedge\limits_{u\in\Phi_L(X)}\mathcal{A}(u)\rightarrow {\rm sub}(u,v).$$

Firstly, for each $a\in L$, by (F2) and (D1), $$w(a_X)=\bigvee\limits_{u\in\Phi_L(X)}\mathcal{A}(u)\wedge u(a_X)\geq\bigvee\limits_{u\in\Phi_L(X)}\mathcal{A}(u)\wedge a=1\wedge a=a.$$For all $A,B\in\mathcal{O}(X)$, clearly we have $w(A\wedge B)\leq w(A)\wedge w(B)$ and
$$\begin{array}{lll}\ &w(A)\wedge w(B)
\\=&\bigvee\limits_{u_1,u_2\in\Phi_L(X)}\mathcal{A}(u_1)\wedge u_1(A)\wedge\mathcal{A}(u_2)\wedge u_2(B)
\\\leq&\bigvee\limits_{u_1,u_2,v\in\Phi_L(X)}\mathcal{A}(v)\wedge{\rm sub}(u_1,v)\wedge {\rm sub}(u_2,v)\wedge u_1(A)\wedge u_2(B)&({\rm D2})
\\\leq&\bigvee\limits_{v\in\Phi_L(X)}\mathcal{A}(v)\wedge v(A)\wedge v(B)\\ =&\bigvee\limits_{v\in\Phi_L(X)}\mathcal{A}(v)\wedge v(A\wedge B)&({\rm F1})\\=&w(A\wedge B).\end{array}$$Hence, $w(A\wedge B)=w(A)\wedge w(B)$. Therefore, $w\in\Phi_L(X)$.

Secondly, for every $v\in\Phi_L(X)$,
$$\begin{array}{llll}{\rm sub}(w,v)&=&\bigwedge\limits_{A\in\mathcal{O}(X)}w(A)\rightarrow v(A)
\\ \ &=&\bigwedge\limits_{A\in\mathcal{O}(X)}\bigwedge\limits_{u\in\Phi_L(X)}(\mathcal{A}(u)\wedge u(A))\rightarrow v(A)&({\rm Prop. 2.2(6)})
\\ \ &=&\bigwedge\limits_{u\in\Phi_L(X)}\mathcal{A}(u)\rightarrow \bigwedge\limits_{A\in\mathcal{O}(X)}u(A)\rightarrow v(A)&({\rm Prop. 2.2(7,12)})
\\ \ &=&\bigwedge\limits_{u\in\Phi_L(X)}\mathcal{A}(u)\rightarrow {\rm sub}(u,v).\end{array}$$

In a summary,  $\sqcup\mathcal{A}=\bigvee\limits_{u\in\Phi_L(X)}\mathcal{A}(u)\wedge u$.}

\begin{Lem}Let $\mathcal{A}$ be a directed $L$-subset of $(\Phi_L(X),{\rm sub})$. Define $\widetilde{\mathcal{A}}:\mathcal{O}(\Phi_L(X))\longrightarrow L$ by $$\widetilde{\mathcal{A}}(W)=\bigvee\limits_{u\in\Phi_L(X)}\mathcal{A}(u)\wedge W(u).$$Then

{\rm(1)} $\widetilde{\mathcal{A}}\in\Phi^2_L(X)$;

{\rm(2)} $\mu_X(\widetilde{\mathcal{A}})=\sqcup\mathcal{A}$.\end{Lem}

\p{(1)
Firstly, for each $a\in L$, by (D1), $$\widetilde{\mathcal{A}}(a_{\Phi_L(X)})=\bigvee\limits_{u\in\Phi_L(X)}\mathcal{A}(u)\wedge a_{\Phi_L(X)}(u)=a\wedge\bigvee\limits_{u\in\Phi_L(X)}\mathcal{A}(u)=a\wedge 1=a.$$
Secondly, for all $W_1,W_2\in\mathcal{O}(\Phi_L(X))$, it is clear that $\widetilde{\mathcal{A}}(W_1\wedge W_2)\leq\widetilde{\mathcal{A}}(W_1)\wedge \widetilde{\mathcal{A}}(W_2)$ and
$$\begin{array}{lll}\ &\widetilde{\mathcal{A}}(W_1)\wedge\widetilde{\mathcal{A}}(W_2)
\\=&\bigvee\limits_{u,v\in\Phi_L(X)}\mathcal{A}(u)\wedge W_1(u)\wedge\mathcal{A}(v)\wedge W_2(v)
\\\leq&\bigvee\limits_{u,v,w\in\Phi_L(X)}\mathcal{A}(w)\wedge{\rm sub}(u,w)\wedge{\rm sub}(v,w)\wedge W_1(u)\wedge W_2(v)&({\rm D2})
\\\leq &\bigvee\limits_{w\in\Phi_L(X)}\mathcal{A}(w)\wedge W_1(w)\wedge W_2(w)
\\=&\bigvee\limits_{w\in\Phi_L(X)}\mathcal{A}(w)\wedge (W_1\wedge W_2)(w)&({\rm F1})
\\=&\widetilde{\mathcal{A}}(W_1\wedge W_2).\end{array}$$ Then, $\widetilde{\mathcal{A}}(W_1)\wedge\widetilde{\mathcal{A}}(W_2)=\widetilde{\mathcal{A}}(W_1\wedge W_2)$. Hence, $\widetilde{\mathcal{A}}\in\Phi_L^2(X)$.

(2) For each $A\in\mathcal{O}(X)$, we have $$\mu_X(\widetilde{\mathcal{A}})(A)=\widetilde{\mathcal{A}}(\phi(A))=\bigvee\limits_{u\in\Phi_L(X)}\mathcal{A}(u)\wedge \phi(A)(u)=\bigvee\limits_{u\in\Phi_L(X)}\mathcal{A}(u)\wedge u(A)=(\sqcup\mathcal{A})(A).$$ Hence, $\mu_X(\widetilde{\mathcal{A}})=\sqcup\mathcal{A}$.}

\begin{Pro}$r:\Phi_L(X)\longrightarrow X$ preserves suprema of directed $L$-subsets.\end{Pro}

\p{Let $\mathcal{A}$ be a directed $L$-subset of $\Phi_L(X)$. Then by Lemma 6.6(2),
$$r(\sqcup\mathcal{A})=r\cdot\mu_X(\widetilde{\mathcal{A}})=r\cdot\Phi_Lr(\widetilde{\mathcal{A}}).$$ We will prove that $r\cdot\Phi_Lr(\widetilde{\mathcal{A}})\leq\sqcup r_L^\rightarrow(\mathcal{A})$ (the inverse inequality is routine).
Since $\sqcup r_L^\rightarrow(\mathcal{A})=\sqcap(r_L^\rightarrow(\mathcal{A}))^u=r([(r_L^\rightarrow(\mathcal{A}))^u])$, we only need to show that $\Phi_Lr(\widetilde{\mathcal{A}})\leq[(r_L^\rightarrow(\mathcal{A}))^u]$.

For every $A\in\mathcal{O}(X)$, $$\Phi_Lr(\widetilde{\mathcal{A}})(A)=\widetilde{\mathcal{A}}(r_L^\leftarrow(A))=\bigvee\limits_{u\in\Phi_L(X)}\mathcal{A}(u)\wedge r_L^\leftarrow(A)(u)=\bigvee\limits_{u\in\Phi_L(X)}\mathcal{A}(u)\wedge A(r(u)).$$
For each $x\in X$,
$$\begin{array}{llll}(r_L^\rightarrow(\mathcal{A}))^u(x)
&=&\bigwedge\limits_{y\in X}r_L^\rightarrow(\mathcal{A})(y)\rightarrow e(y,x)
\\&=&\bigwedge\limits_{u\in\Phi_L(X)}\mathcal{A}(u)\rightarrow e(r(u),x)
\\&=&\bigwedge\limits_{u\in\Phi_L(X)}\mathcal{A}(u)\rightarrow\bigwedge\limits_{B\in\mathcal{O}(X)}B(r(u))\rightarrow B(x)
\\&=&\bigwedge\limits_{u\in\Phi_L(X)}\bigwedge\limits_{B\in\mathcal{O}(X)}(\mathcal{A}(u)\wedge B(r(u)))\rightarrow B(x)&{\rm (Prop.2.2(7))}
\\&=&\bigwedge\limits_{B\in\mathcal{O}(X)}(\bigvee\limits_{u\in\Phi_L(X)}\mathcal{A}(u)\wedge B(r(u)))\rightarrow B(x).&{\rm (Prop.2.2(6))}\end{array}$$
and $$\begin{array}{lll}\ &[(r_L^\rightarrow(\mathcal{A}))^u](A)={\rm sub}((r_L^\rightarrow(\mathcal{A}))^u,A)
\\=&\bigwedge\limits_{x\in X}(r_L^\rightarrow(\mathcal{A}))^u(x)\rightarrow A(x)
\\=&\bigwedge\limits_{x\in X}\{\bigwedge\limits_{B\in\mathcal{O}(X)}(\bigvee\limits_{u\in\Phi_L(X)}\mathcal{A}(u)\wedge B(r(u)))\rightarrow B(x)\}\rightarrow A(x)
\\\geq&\bigwedge\limits_{x\in X}\{(\bigvee\limits_{u\in\Phi_L(X)}\mathcal{A}(u)\wedge A(r(u)))\rightarrow A(x)\}\rightarrow A(x)&{\rm (Prop.2.2(6))}
\\\geq&\bigvee\limits_{u\in\Phi_L(X)}\mathcal{A}(u)\wedge A(r(u))&{\rm (Prop.2.2(4))}
\\=&\Phi_Lr(\widetilde{\mathcal{A}})(A).\end{array}$$
This completes the proof.}

For every $u\in\Phi_L(X)$, define $\mathcal{A}_u\in L^{\Phi_L(X)}$ by
$$\mathcal{A}_u(v)
=\bigvee\limits_{A\in\mathcal{O}(X)}u(A)\wedge {\rm sub}(v,[A])\ (\forall v\in\Phi_L(X)).$$

\begin{Lem}$\mathcal{A}_u$ is directed in $(\Phi_L(X),{\rm sub})$ and $\sqcup\mathcal{A}_u=u$.\end{Lem}

\p{(1) Firstly, by (F2), $$\bigvee\limits_{v\in\Phi_L(X)}\mathcal{A}_u(v)
=\bigvee\limits_{v\in\Phi_L(X)}\bigvee\limits_{A\in\mathcal{O}(X)}u(A)\wedge {\rm sub}(v,[A])\geq u(1_X)\wedge{\rm sub}([1_X],[1_X]))
=1.$$

Secondly, for all $v_1,v_2\in\Phi_L(X)$,
$$\begin{array}{lll}\ &\mathcal{A}_u(v_1)\wedge\mathcal{A}_u(v_2)
\\=&\bigvee\limits_{A_1,A_2\in\mathcal{O}(X)}u(A_1)\wedge {\rm sub}(v_1,[A_1])\wedge u(A_2)\wedge {\rm sub}(v_2,[A_2])
\\=& \bigvee\limits_{A_1,A_2\in\mathcal{O}(X)}u(A_1\wedge A_2)\wedge {\rm sub}(v_1,[A_1])\wedge {\rm sub}(v_2,[A_2])&({\rm F1})
\\\leq& \bigvee\limits_{A_1,A_2\in\mathcal{O}(X)}u(A_1\wedge A_2)\wedge {\rm sub}(v_1,[A_1\wedge A_2])\wedge {\rm sub}(v_2,[A_1\wedge A_2])&({\rm Prop. 2.14})
\\\leq& \bigvee\limits_{A\in\mathcal{O}(X)}u(A)\wedge {\rm sub}(v_1,[A])\wedge {\rm sub}(v_2,[A])\\=&\bigvee\limits_{A\in\mathcal{O}(X)}u(A)\wedge {\rm sub}([A],[A])\wedge {\rm sub}(v_1,[A])\wedge {\rm sub}(v_2,[A])
\\\leq&\bigvee\limits_{v\in\Phi_L(X)}\bigvee\limits_{A\in\mathcal{O}(X)}u(A)\wedge {\rm sub}(v,[A])\wedge {\rm sub}(v_1,v)\wedge {\rm sub}(v_2,v)
\\=&\bigvee\limits_{v\in\Phi_L(X)}\mathcal{A}_u(v)\wedge {\rm sub}(v_1,v)\wedge {\rm sub}(v_2,v).\end{array}$$

Hence $\mathcal{A}_u$ is directed.

Thirdly, by Lemma 6.5, for every $B\in\mathcal{O}(X)$, we have $$\sqcup\mathcal{A}_u(B)=\bigvee\limits_{v\in\Phi_L(X)}\mathcal{A}_u(v)\wedge v(B)=\bigvee\limits_{v\in\Phi_L(X)}\bigvee\limits_{A\in\mathcal{O}(X)}u(A)\wedge {\rm sub}(v,[A])\wedge v(B).$$
By (F3$'$), we only need to show that $\bigvee\limits_{v\in\Phi_L(X)}{\rm sub}(v,[A])\wedge v(B)=[A](B)$.

Firstly, by Proposition 2.2(3), $$\bigvee\limits_{v\in\Phi_L(X)}{\rm sub}(v,[A])\wedge v(B)\leq\bigvee\limits_{v\in\Phi_L(X)}(v(B)\rightarrow[A](B))\wedge v(B)\leq [A](B).$$
Secondly, $$\bigvee\limits_{v\in\Phi_L(X)}{\rm sub}(v,[A])\wedge v(B)\geq{\rm sub}([A],[A])\wedge [A](B)=[A](B)$$
Thus, $\bigvee\limits_{v\in\Phi_L(X)}{\rm sub}(v,[A])\wedge v(B)=[A](B)$.

Hence, $\sqcup\mathcal{A}_u=u$.}

\begin{Pro}For $u\in\Phi_L(X)$, $r(u)=\sqcup u^l$.\end{Pro}

\p{We need to show that $e(r(u),x)={\rm sub}(u^l,{\downarrow}x)$ for all $x\in X$.

First of all, $$\begin{array}{lll}\ &{\rm sub}(u^l,{\downarrow}x)\\=&\bigwedge\limits_{y\in X}u^l(y)\rightarrow e(y,x)
\\=&\bigwedge\limits_{y\in X}\bigwedge\limits_{A\in\mathcal{O}(X)}(u(A)\wedge {\rm sub}(A,{\uparrow}y))\rightarrow e(y,x)&({\rm Prop. 2.2(6)})
\\=&\bigwedge\limits_{y\in X}\bigwedge\limits_{A\in\mathcal{O}(X)}(u(A)\wedge {\rm sub}(y,\sqcap A))\rightarrow e(y,x)&({\rm Rem.2.10(2)})
\\=&\bigwedge\limits_{A\in\mathcal{O}(X)}u(A)\rightarrow (\bigwedge\limits_{y\in X}e(y,\sqcap A)\rightarrow e(y,x))&({\rm Prop. 2.2(7,12)})
\\=&\bigwedge\limits_{A\in\mathcal{O}(X)}u(A)\rightarrow (e(y,\sqcap A)\rightarrow e(y,x))&({\rm Prop. 2.2(7,12)})
\\=&\bigwedge\limits_{A\in\mathcal{O}(X)}u(A)\rightarrow e(\sqcap A,x).&({\rm Prop. 2.7(3)})\end{array}$$
Then we shall prove that $e(r(u),x)=\bigwedge\limits_{A\in\mathcal{O}(X)}u(A)\rightarrow e(\sqcap A,x)$ for all $x\in X$.

Firstly, for all $A\in\mathcal{O}(X)$, $$\begin{array}{llll}e(r(u),x)\wedge u(A)&=&e(r(u),x)\wedge {\rm sub}([A],u)&({\rm F4}')\\\ &\leq& e(r(u),x)\wedge e(r([A]),r(u))&({\rm Lam. 6.2(2)})\\\ &=&e(r(u),x)\wedge e(\sqcap A,r(u))&({\rm Prop. 6.4})\\\ &\leq& e(\sqcap A,x).&({\rm E3})\end{array}$$ By the arbitrariness of $A\in\mathcal{O}(X)$, we have $$e(r(u),x)\leq\bigwedge\limits_{A\in\mathcal{O}(X)}u(A)\rightarrow e(\sqcap A,x).$$

Secondly,  for every $x\in X$,
$$\begin{array}{lll}
\ &e(r(u),x)\\=&e(r(\sqcup\mathcal{A}_u),x)=e(\sqcup r_L^\rightarrow(\mathcal{A}_u),x)&({\rm Prop. 6.7,\ Lem. 6.8})
\\=&\bigwedge\limits_{y\in X}r_L^\rightarrow(\mathcal{A}_u)(y)\rightarrow e(y,x)&({\rm Rem.2.10(1)})
\\=&\bigwedge\limits_{y\in X}\bigwedge\limits_{r(v)=y}\mathcal{A}_u(v)\rightarrow e(y,x)&({\rm Prop. 2.2(6), Rem.2.10(1)})
\\=&\bigwedge\limits_{v\in \Phi_L(X)}\bigwedge\limits_{A\in O(X)}(u(A)\wedge {\rm sub}(v,[A]))\rightarrow e(r(v),x)&({\rm Prop. 2.2(7)})
\\=&\bigwedge\limits_{A\in O(X)}u(A)\rightarrow (\bigwedge\limits_{v\in\Phi_L(X)}{\rm sub}(v,[A])\rightarrow e(r(v),x))&({\rm Prop. 2.2(7,12)}).\end{array}$$
We only need to show that $e(\sqcap A,x)\wedge {\rm sub}(v,[A])\leq e(r(v),x)$ for all $v\in \Phi_L(X)$. In fact, by Lemma 6.2(2) and (E3), $$\begin{array}{lll}e(\sqcap A,x)\wedge {\rm sub}(v,[A])&\leq& e(\sqcap A,x)\wedge e(r(v),r([A]))\\ \ &=&e(\sqcap A,x)\wedge e(r(v),\sqcap A)\leq e(r(v),x).\end{array}$$ This completes the proof.}

\begin{Pro}$(X,e_{\mathcal{O}(X)})$ is an $L$-valued continuous lattice.\end{Pro}

\p{By Proposition 6.9, $x=r([x])=\sqcup[x]^l$ holds for every $x\in X$. By Propositions 6.4 and 4.10, $X$ is an $L$-valued continuous lattice.}

\end{document}